

%



\input mssymb

\iftrue
\else
\baselineskip=1.67\normalbaselineskip     
\smallskipamount=1.33\smallskipamount     
\medskipamount=1.33\medskipamount         
\bigskipamount=1.33\bigskipamount
\magnification=\magstep1

\hsize=5.9truein
\hoffset=0.5truein
\catcode`@=11
\def\plainoutput{\shipout\vbox{\makeheadline\pagebody\vskip10pt\makefootline}%
  \advancepageno
  \ifnum\outputpenalty>-\@MM \else\dosupereject\fi}
\catcode`@=12
\fi

\iftrue
\magnification=\magstep1
\hsize=6.5truein
\hoffset=0.0truein
\baselineskip 1.4\normalbaselineskip

\tolerance=10000
\def\sqr{$\vcenter{\hrule height .3mm
\hbox {\vrule width .3mm height 2mm \kern 2mm
\vrule width .3mm} \hrule height .3mm}$}
\else

\baselineskip=1.33\normalbaselineskip
\tolerance=10000
\def\sqr{\vcenter{\hrule height .3mm
\hbox {\vrule width .3mm height 2mm \kern 2mm
\vrule width .3mm} \hrule height .3mm}}
\fi

\def \title{\medskip\centerline}
\def \Proof{\noindent {\bf Proof.\ \ }}
\def \R{{\Bbb R}}
\def \N{{\Bbb N}}

\def\CC{{\Bbb C}}
\def \nm#1{\left\|#1\right\|}
\def \bgnm#1{\biggl\|#1\biggr\|}
\def \nmm#1{\left|\!\left|\!\left|#1\right|\!\right|\!\right|}
\def \normm#1{|\!|\!|#1|\!|\!|}
\def \nmm{\normm}
\def \seq#1#2{#1_1,\dots,#1_#2}
\def \sleq#1#2{#1_1<\dots<#1_#2}
\def \speq#1#2{#1_1+\dots+#1_#2}
\def \sm#1#2{\sum_{#1=1}^#2}

\def \G{\Gamma}
\def \d{\delta}
\def \e{\epsilon}
\def\pl{\char'40l} 
\def \l{\lambda}

\def \s{\sigma}

\def\bfe{{\bf e}}

\def \ran{{\rm ran}}
\def \supp{{\rm supp}}
\def \exp{{\rm exp}}
\def \supp{{\rm supp}}
\def \ra{\rightarrow}
\def \sp#1{\langle#1\rangle}

\font\cbf=cmssbx10
\font\ncbf=cmssbx10 at 9pt 
\def\abstract#1{{\ninepoint{\narrower\smallskip\noindent
	{\ncbf Abstract.} #1\smallskip}\vskip.5truein}}  
	
\def\hangbox to #1 #2{\vskip1pt\hangindent #1\noindent \hbox to #1{#2}$\!\!$}
\def\myitem#1{\hangbox to 40pt {#1\hfill}} 
\newskip\ttglue 
\def\ninepoint{\def\rm{\fam0\ninerm}
  	\textfont0=\ninerm \scriptfont0=\sixrm \scriptscriptfont0=\fiverm
  	\textfont1=\ninei  \scriptfont1=\sixi  \scriptscriptfont1=\fivei
  	\textfont2=\ninesy  \scriptfont2=\sixsy  \scriptscriptfont2=\fivesy
	\textfont3=\tenex  \scriptfont3=\tenex  \scriptscriptfont3=\tenex
	\textfont\itfam=\nineit  \def\it{\fam\itfam\nineit}
	\textfont\slfam=\ninesl  \def\sl{\fam\slfam\ninesl}
	\textfont\ttfam=\ninett  \def\tt{\fam\ttfam\ninett}
	\textfont\bffam=\ninebf  \scriptfont\bffam=\sixbf
	\scriptscriptfont\bffam=\fivebf  \def\bf{\fam\bffam\ninebf}
	\tt  \ttglue=.5em plus.25em minus.15em
	\normalbaselineskip=11pt
	\setbox\strutbox=\hbox{\vrule height8pt depth3pt width0pt}
	\let\sc=\sevenrm  \let\big=\ninebig \normalbaselines\rm}

	\font\ninerm=cmr9 \font\sixrm=cmr6 \font\fiverm=cmr5
	\font\ninei=cmmi9  \font\sixi=cmmi6   \font\fivei=cmmi5
	\font\ninesy=cmsy9  \font\sixsy=cmsy6 \font\fivesy=cmsy5
	\font\nineit=cmti9  \font\ninesl=cmsl9  \font\ninett=cmtt9
	\font\ninebf=cmbx9  \font\sixbf=cmbx6 \font\fivebf=cmbx5
	\def\ninebig#1{{\hbox{$\textfont0=\tenrm\textfont2=\tensy
	\left#1\vbox to7.25pt{}\right.$}}}

\title {\bf THE UNCONDITIONAL BASIC SEQUENCE PROBLEM}
\smallskip
\centerline {W. T. Gowers and B. Maurey - preliminary version}
\bigskip

\abstract{We construct a Banach space that does not contain
any infinite unconditional basic sequence.}

\noindent {\bf \S 0. Introduction.}

A fundamental role in the theory of Banach spaces is played by the notion of a
Schauder basis. If $X$ is a Banach space, then a sequence $(x_n)_{n=1}^\infty$
is a Schauder basis (or simply a basis) of $X$ if the closed linear span of
$\{x_n\}_{n=1}^\infty$ is $X$ and $\sum_{n=1}^\infty a_nx_n$ is zero only if
each $a_n$ is zero. The second condition, asserting a kind of independence,
clearly depends very much on the order of the $x_n$, and it is certainly
possible for a permutation of a basis to fail to be a basis. On the other hand,
many bases that occur naturally, such as the standard basis of $\ell_p$ when
$1\le p<\infty$, are bases under any permutation. It is therefore natural to
give a name to this special kind of basis. As it happens there are several
equivalent definitions.

\noindent {\bf Theorem 0$\cdot$1.}  {\sl Let $X$ be a (real or complex) Banach
space and let $(x_n)_{n=1}^\infty$ be a basis of~$X$. Then the following are
equivalent.

\item {(i)} $(x_{\pi(n)})_{n=1}^\infty$ is a basis of $X$ for every permutation
$\pi:\N\ra\N$.
\item {(ii)} Sums of the form $\sum_{n=1}^\infty a_nx_n$ converge
unconditionally whenever they converge.
\item {(iii)} There exists a constant $C$ such that, for every sequence of
scalars $(a_n)_{n=1}^\infty$ and every sequence of scalars
$(\e_n)_{n=1}^\infty$ of modulus at most $1$, we have the inequality
$$\bgnm{\sum_{n=1}^\infty\e_na_nx_n}\le C\bgnm{\sum_{n=1}^\infty a_nx_n}\ .$$}

A basis satisfying these conditions is called an {\sl unconditional} basis, and
a basis satisfying the third condition for some given constant $C$ is called
$C$-{\sl unconditional}. An infinite sequence that is a basis of its closed
linear span is called a basic sequence: if it is an unconditional basis of its
closed linear span then it is an unconditional basic sequence.

For a long time a major unsolved problem was whether every separable Banach
space had a basis. This question was answered negatively by
P.~Enflo in 1973 [E].
On the other hand, it is not hard to show that every space contains a basic
sequence. Spaces with {\it unconditional} bases have much more structure than
general spaces, so examples of spaces without them are easier to find. Indeed,
the spaces $C([0,1])$ and  $L_1$
do not have an unconditional basis. This leaves the question
of whether every space contains an unconditional basic sequence, or
equivalently has an infinite-dimensional subspace with an unconditional basis.
The earliest reference we know for the problem is 
Bessaga-Pe\pl czy\'nski (1958), where it appears as problem 5.1;
actually we solve the more precise problem 5.11, since our example
is a reflexive Banach space. The easier related problem 5.12
was solved already many years ago [MR].

In the summer of 1991, the first named author found a counterexample. A short
time afterwards the second named author independently found a counterexample as
well. On comparing our examples, we discovered that they were almost identical,
as were the proofs that they were indeed counterexamples, so we decided to
publish jointly and work together on further properties of the space. As a
result of our collaboration, the proofs of some of the main lemmas have been
simplified and tightened.

After reading our original preprints, W. B. Johnson pointed out that our
proof(s) could be modified to give a much stronger property of the space. J.
Lindenstrauss had asked whether every infinite-dimensional Banach space was
{\sl decomposable}, that is, could be written as a topological direct sum
$Y\oplus Z$ with $Y$ and $Z$ infinite-dimensional. Johnson's observation was
that our space, which for the remainder of the introduction we shall call $X$,
is not only not decomposable, but does not even have a decomposable subspace.
That is, $X$ is {\sl hereditarily indecomposable} or H.I. Equivalently, if $Y$
and $Z$ are two infinite-dimensional subspaces of $X$ and $\e>0$ then there
exist $y\in Y$ and $z\in Z$ such that $\nm y=\nm z=1$ and $\nm{y-z}<\e$. This
turned out to be a key property of $X$ in that all of the pathological
properties that we know about $X$ can be deduced from the fact that it is H.I.
In particular it is easy to see that a space with this property cannot contain
an unconditional basic sequence.

Another immediate consequence is that either the space is a new prime Banach
space (which means that it is isomorphic to all its complemented subspaces) or
it fails to be isomorphic to a subspace of finite codimension. If the second
statement is true then it must fail to be isomorphic to a subspace of
codimension 1, giving a counterexample to a question of Banach which has come
to be known as the hyperplane problem. The first author modified the
construction of $X$ to give such a counterexample, and in fact one with an
unconditional basis [G]. Soon afterwards, we managed to use the H.I. property
to show that the complex version of $X$ gives another example. Later, we were
able to pass to the real case, so $X$ itself is a counterexample to the
hyperplane problem by virtue of being a H.I. space.

In fact, the space of operators on~$X$ is very small: every bounded linear
operator on~$X$ can be written as $\l Id+S$, where $S$ is a strictly singular
operator. A question we have not answered is whether there exists a space on
which every bounded linear operator is of the form $\l Id +K$ for a compact
operator $K$. We do not even know whether our space has that property, though
it seems unlikely.

The rest of this paper is divided into five sections. The first concerns the
notion of an asymptotic set, which is a definition of great importance for this
problem, but which arises most naturally in the context of the distortion
problem, about which we shall have more to say later. In particular, we give a
criterion for a space to have an equivalent norm in which it contains no
$C$-unconditional basic sequence.

The second section is about a remarkable space constructed by T. Schlumprecht,
on which our example builds. We show that, for every $C$, his space satisfies
our criterion and therefore can be renormed so as not to contain a
$C$-unconditional basic sequence.

The third section contains the definition of $X$ and a proof that it is H.I.
and therefore contains no unconditional basic sequence,
and ends with the (easy) proof that $X$ is reflexive. The fourth is about
consequences of this property, especially the existence of very few operators
on a complex space having it. The final section concerns the passage to the
real case of the results of the previous one.

We are very grateful to W. B. Johnson for his influence on this paper. As we
have mentioned, his observation that our space is H.I. lies at the heart of all
its interesting properties. He also explained to us much simpler arguments for
some of the consequences of this property. We would also like to thank P. G.
Casazza and R. G. Haydon for interesting conversations and suggestions about
the problems solved here.

For the rest of this paper we shall use the words ``space'' and ``subspace'' to
refer to infinite-dimensional spaces and subspaces. Similarly a basis will
always be assumed to be infinite.
\bigbreak

\noindent {\bf \S 1. Asymptotic sets.}

Let $X$ be a normed space and let $S(X)$ be its unit sphere. We shall say that
a subset $A\subset S(X)$ is {\sl asymptotic} if $A\cap S(Y)\ne\emptyset$ for
every infinite-dimensional (not necessarily closed) subspace $Y\subset X$. A
key observation for this paper is that if a space $X$ contains infinitely many
asymptotic sets that are all well disjoint from one another, then these can be
used to construct an equivalent norm on $X$ such that no sequence is
$C$-unconditional in this norm. In this section, we shall make that statement
precise and prove it. The technique we use will underlie our later arguments as
well.

Let $A_1,A_2,\dots$ be a sequence of subsets of the unit sphere of a normed
space $X$ and let $A_1^*,A_2^*,\dots$ be a sequence of subsets of the unit ball
of $X^*$. (It is slightly more convenient in applications to take the ball
rather than the sphere). We shall say that $A_1,A_2,\dots$ and
$A_1^*,A_2^*,\dots$ are an {\sl asymptotic biorthogonal system with constant}
$\d$ if the following conditions hold.

\item {(i)} For every $n\in\N$, the set $A_n$ is asymptotic.

\item {(ii)} For every $n\in\N$ and every $x\in A_n$ there exists $x^*\in
A_n^*$ such that $x^*(x)>1-\d$.

\item {(iii)} For every $n,m\in\N$ with $n\ne m$, every $x\in A_n$ and every
$x^*\in A_m^*$, $|x^*(x)|<\d$.

Under these circumstances, we shall say that $X$ contains an asymptotic
biorthogonal system. The definition is not interesting if $\d>1/2$ since one
may take $A_n=S(X)$ and $A_n^*={1\over 2}B(X^*)$ for every $n$. On the other
hand, if $\d<1/2$, it is not at all obvious that any Banach space contains an
asymptotic biorthogonal system with constant $\d$. We shall see later, however,
that this is not as rare a phenomenon as it might seem.

Note that the $A_n$ are separated in the following sense. If $n\ne m$, $x\in
A_n$ and $y\in A_m$, then there exists $x^*\in A_n^*$ such that $x^*(x)\ge
1-\d$ and $|x^*(y)| <\d$. Since $A_n^*\subset B(X^*)$, it follows that
$\nm{x-y}\ge 1-2\d$.

The main result of this section is the following theorem.

\proclaim Theorem 1$\cdot$1. Let $0<\d<1/36$ and let $X$ be a separable normed
space containing an asymptotic biorthogonal system with constant $\d$. Then
there is an equivalent norm on $X$ such that no sequence is
$1/\sqrt{36\d}$-unconditional.

\Proof   Let $\nm.$ be the original norm on $X$ and let $A_1,A_2,\dots$ and
$A_1^*,A_2^*,\dots$ be the asymptotic biorthogonal system with constant $\d$.
For each $n$ let $Z_n^*$ be a countable subset of $A_n^*$ such that
property (ii) of asymptotic biorthogonal systems holds for some $x^*\in Z_n^*$.
Let $Z^*=\bigcup_{n=1}^\infty Z_n^*$. Next, let $\s$ be an injection to the
natural numbers from the collection of finite sequences of elements of $Z^*$.

We shall now define a collection of functionals which we call {\sl special}
functionals. A {\sl special sequence} of functionals {\sl of length $r$} is a
sequence of the form $z_1^*,z_2^*,\dots,z_r^*$, where $z_1^*\in Z_1^*$ and, for
$1\le i<r$, $z_{i+1}^*\in Z_{\s(z_1^*,\dots,z_i^*)}^*$. A {\sl special
functional of length $r$} is simply the sum of a special sequence of length
$r$. We shall let $\G_r$ stand for the collection of special functionals of
length $r$.

We can now define an equivalent norm on $X$. Let $r=\lfloor\d^{-1/2}\rfloor$
and define a norm $\nmm .$ by
$$\nmm x=\nm x\vee r\max\Bigl\{|z^*(x)|:z^*\in \G_r\Bigr\}\ .$$

Let $x_1,x_2,\dots$ be any sequence of linearly independent vectors in $X$. We
shall show that it is not
$(r-1)/4$-unconditional in the norm $\nmm.$. We shall do
this by constructing a sequence of vectors $\seq z r$, generated by
$x_1,x_2,\dots$ and disjointly supported with respect to these vectors, with
the property that
$$(r-1)\nmm{\sum_{i=1}^r(-1)^iz_i} < 4 \nmm{\sum_{i=1}^rz_i}\ .$$
This will obviously prove the theorem, since
$(r-1)/4 > 1/\sqrt{36 \delta}$.

Let $X_1$ be the algebraic subspace generated by $(x_i)_{i=1}^\infty$. Since
$A_1$ is an asymptotic set, we can find $z_1\in A_1\cap X_1$. This implies that
$z_1$ has norm 1 and is generated by finitely many of the $x_i$. Next we can
find $z_1^*\in Z_1^*$ such that $z_1^*(z_1)>1-\d$. Now let $X_2$ be the
algebraic subspace generated by all the $x_i$ not used to generate $z_1$. Since
$A_{\sigma(z_1^*)}$ is asymptotic, we can find $z_2\in A_{\s(z_1^*)}\cap X_2$
of norm 1. We can then find $z_2^*\in Z_{\s(z_1^*)}^*$ such that
$z_2^*(z_2)>1-\d$.

Continuing this process, we obtain sequences $\seq z r$ and $z_1^*,\dots,z_r^*$
with the following properties. First, $\nm{z_i}=1$ for each $i$. Second,
$z_{i+1}^*\in Z_{\s(z_1^*,\dots,z_i^*)}^*$
for each $i$ (i.e. $z_1^*,\dots,z_r^*$ is a special sequence of length $r$).
Third, $z_i^*(z_i)>1-\d$ for each $i$. Fourth, since $\s$ is an injection, the
$z_i^*$ belong to different $A_n^*$s, so $|z_i^*(z_j)|<\d$ when $i\ne j$.

Let us now estimate the norm of $\sum_{i=1}^rz_i$. Since $z_1^*,\dots,z_r^*$ is
a special functional of length $r$, the norm is at least
$$\eqalign{r\Bigl(\sum_{i=1}^rz_i^*\Bigr)\Bigl(\sum_{i=1}^rz_i\Bigr)
& >  r\bigl((1-\d)r-\d r(r-1)\bigr) \cr
&\ge r(r-1)\ .\cr}$$

On the other hand, if $(w_i^*)_{i=1}^r$ is any special sequence of length $r$,
let $t$ be maximal such that $w_i^*=z_i^*$ (or zero if $w_1^*\ne z_1^*$). Then
$$\Bigl|\sum_{i=1}^r(-1)^iw_i^*(z_i)\Bigr|\le\Bigl|\sum_{i=1}^t(-1)^iw_i^*(z_i)
\Bigr| +|w_{t+1}^*(z_{t+1})|+\sum_{i=t+2}^r|w_i^*(z_i)|\ .$$
Since $\s$ is an injection, $w_i^*$ and $z_j^*$ are chosen from different sets
$A_n^*$ whenever $i\ne j$ or $i=j>t+1$. By property (iii) this tells us that
$|w_i^*(z_j)|<\d$.
In particular, $\sum_{i=t+2}^r|w_i^*(z_i)|<\d r$. When $i<t$
we know that $1-\d<w_i^*(z_i)\le 1$, so
$\Bigl|\sum_{i=1}^t(-1)^iw_i^*(z_i)\Bigr|\le 1+\d t/2$. It follows that
$$\Bigl|\sum_{i=1}^r(-1)^iw_i^*(z_i)\Bigr|\le 1+\d r/2+1+\d r\le 2(1+\d r)\ .$$
We also know that $\sum_{i\ne j}|w_i^*(z_j)|\le\d r(r-1)$. Finally, by the
triangle inequality, $\nm{\sum_{i=1}^r(-1)^iz_i}\le r$.

Putting all these estimates together, we find that
$$\nmm{\sum_{i=1}^r(-1)^iz_i}\le r\Bigl(2(1+\d r)+\d r(r-1)\Bigr) < 4r$$
from which it follows that the basic sequence $x_1,x_2,\dots$ was not
$(r-1)/4$-unconditional in the equivalent norm.
\hfill $\square$
\bigskip

With a little more care one can increase the best unconditional constant from
roughly $\d^{-1/2}$ to roughly $\d^{-1}$, but some of the details of this would
obscure the main point of the proof. It also does not seem to be necessary in
applications. Indeed, it is not known whether there exists a space containing
an asymptotic biorthogonal system for some (non-trivial) $\d$ but not for every
$\d>0$. In the next section, we examine a space that contains them for every
$\d$.
\bigbreak

\noindent {\bf \S 2. Schlumprecht's space.}

A space $(Y,\nm.)$ is said to be $\l$-{\sl distortable} if there exists an
equivalent norm $\nmm.$ on $Y$ such that, for every subspace $Z\subset Y$ the
quantity $\sup\{\nmm y/\nmm z:\nm y=\nm z=1\}$ is at least~$\l$. A space is
{\sl distortable} if it is $\l$-distortable for some $\l>1$. A famous open
problem, known as the distortion problem, used to be whether $\ell_2$ was
distortable. This is equivalent to asking whether every space isomorphic to
$\ell_2$ contains a subspace almost isometric to~$\ell_2$. A few months after
the results in this paper were obtained, the distortion problem was also
solved, by E. Odell and T. Schlumprecht [OS]; actually a stronger statement is
proved in [OS], namely that $\ell_2$ contains an asymptotic biorthogonal system
with any constant $\d>0$. This implies that $\ell_2$ can be renormed so as not
to contain a $C$-unconditional basic sequence. However, we shall consider in
this section a space constructed by Schlumprecht~[S1]. This space was the first
known example of a space that is $\l$-distortable for every $\l$. The main
result of this section is that it contains an asymptotic biorthogonal system
for any $\d$. In proving this, we shall use very little more than what was
already proved by Schlumprecht in order to show that it is arbitrarily
distortable.

First, let us give the definition of Schlumprecht's space. He defines a class
of functions $f:[1,\infty)\ra [1,\infty)$, which we shall call $\cal F$, as
follows. The function $f$ is a member of $\cal F$ if it satisfies the following
five conditions.

\item {(i)} $f(1)=1$ and $f(x)<x$ for every $x>1$;
\item {(ii)} $f$ is strictly increasing and tends to infinity;
\item {(iii)} $\lim_{x\ra\infty}x^{-q}f(x)=0$ for every $q>0$;
\item {(iv)} the function $x/f(x)$ is concave and non-decreasing;
\item {(v)} $f(xy)\le f(x)f(y)$ for every $x,y\ge 1$.

It is easily verified that $f(x)=\log_2(x+1)$ satisfies these conditions, as
does the function $\sqrt{f(x)}$.

Schlumprecht's space is a Tsirelson-type construction, in that it is defined
inductively. As with an earlier construction due to L. Tzafriri, the
admissibility condition used in Tsirelson's space ([T]) is not needed
(see [CS]). Before giving
the definition, let us fix some notation.

Let $c_{00}$ be the space of sequences of real numbers all but finitely many of
which are zero. We shall let $\bfe_1,\bfe_2,\dots$ stand for the unit vector
basis of this vector space. If $E\subset\N$, then we shall also use the letter
$E$ for the projection from $c_{00}$ to $c_{00}$ defined by
$E\Bigl(\sum_{i=1}^\infty a_i\bfe_i\Bigr)=\sum_{i\in E}a_i\bfe_i$. If
$E,F\subset\N$, then we write $E<F$ to mean that $\max E<\min F$, and if
$k\in\N$ and $E\subset\N$, then we write $k<E$ to mean $k<\min E$ . The {\sl
support} of a vector $x=\sum_{i=1}^\infty x_i\bfe_i\in c_{00}$ is the set of
$i\in\N$ for which $x_i\ne 0$. An {\sl interval} of integers is a subset of
$\N$ of the form $\{a,a+1,\dots,b\}$ for some $a,b\in\N$. We shall also define
the {\sl range} of a vector, written $\ran(x)$, to be the smallest interval
containing its support. We shall write $x<y$ to mean $\ran(x)<\ran(y)$. If
$\sleq x n$ we shall say that $\seq x n$ are {\sl successive}.

Now let $f(x)$ be the function $\log_2(x+1)$ as above. If $x\in c_{00}$, its
norm in Schlumprecht's space is defined inductively by $$\nm x=\nm
x_\infty\vee\sup f(N)^{-1}\sum_{i=2}^N\nm{E_ix}$$ where the supremum runs over
all integers $N\ge 2$ and all sequences of sets $\sleq E N$. Note that this
definition, although apparently circular, in fact determines a unique norm.
Note also that the standard basis of $c_{00}$ is 1-unconditional in this norm,
so there is no difference if we assume that all the sequences $\sleq E N$ are
sequences of {\it intervals}. Later in the paper it will make a great
difference, and we now adopt the convention that all such sequences mentioned
are sequences of intervals.

We now prove various lemmas about this space. As we have already said, they are
essentially due to Schlumprecht [S1][S2],
but are stated here in slightly greater
generality so that they can be applied in the main space of this paper.

Let $\cal X$ be the set of normed spaces of the form $X=(c_{00},\nm.)$ such that
$(\bfe_i)_{i=1}^\infty$ is a normalized monotone basis of $X$. If $f\in\cal F$,
$X\in\cal X$ and every $x\in X$ satisfies the inequality
$$\nm x\ge\sup\biggl\{f(N)^{-1}\sum_{i=1}^N\nm{E_ix}:N\in\N,\sleq E N\biggr\}$$
then we shall say that $X$ {\sl satisfies a lower} $f$-{\sl estimate}. (It is
important that, in the supremum above, the $E_i$ are intervals). Note that this
implies that $\nm{Ex}\le\nm x$ for every interval $E$ and vector $x$, so the
standard basis of a space with a lower $f$-estimate is automatically
bimonotone.

Given a space $X\in\cal X$ and a vector $x\in X$, we shall say that $x$ is an
$\ell_{1+}^n$-{\sl average with constant} $C$ if $\nm x=1$ and
$x=\sum_{i=1}^nx_i$ for some sequence $\sleq x n$ of non-zero elements of $X$
such that $\nm{x_i}\le Cn^{-1}$ for each $i$. An $\ell_{1+}^n$-vector is any
positive multiple of an $\ell_{1+}^n$-average. In other words, a vector $x$ is
an $\ell_{1+}^n$-vector with constant $C$ if it can be written $x=\speq x n$,
where $\sleq x n$, the $x_i$ are non-zero and $\nm{x_i}\le Cn^{-1}\nm x$ for
every $i$.

Finally, by a {\sl block basis} in a space $X\in\cal X$ we mean a sequence
$x_1,x_2,\dots$ of successive non-zero vectors in $X$ (note that such a
sequence must be a basic sequence) and by a {\sl block subspace} of a space
$X\in\cal X$ we mean a subspace generated by a block basis.

\proclaim Lemma 1. Let $f\in\cal F$ and let $X\in\cal X$ satisfy a lower
$f$-estimate. Then, for every $n\in\N$ and every $C>1$, every block subspace
$Y$ of $X$ contains an $\ell_{1+}^n$-average with constant $C$.

\Proof  Suppose the result is false. Let $k$ be an integer such that $k\log
C>\log f(n^k)$ (such an integer exists because of property (iii) in the
definition of $\cal F$), let $N=n^k$, let $\sleq x N$ be any sequence of
successive norm-1 vectors in $Y$ and let $x=\sum_{i=1}^Nx_i$. For every $0\le
i\le k$ and every $1\le j\le n^{k-i}$, let
$x(i,j)=\sum_{t=(j-1)n^i+1}^{jn^i}x_t$. Thus $x(0,j)=x_j$, $x(k,1)=x$ and, for
$1\le i\le k$, each $x(i,j)$ is a sum of $n$ successive $x(i-1,j)$s. By our
assumption, no $x(i,j)$ is an $\ell_{1+}^n$-vector with constant $C$. It
follows easily by induction that $\nm{x(i,j)}\le C^{-i}n^i$, and, in
particular, that $\nm x\le C^{-k}n^k=C^{-k}N$. However, it follows from the
fact that $X$ satisfies a lower $f$-estimate that $\nm x\ge Nf(N)^{-1}$. This
is a contradiction, by choice of $k$. \hfill $\square$ \bigskip

\proclaim Lemma 2. Let $M,N\in\N$ and $C\ge 1$, let $x$ be an
$\ell_{1+}^N$-vector with constant $C$ and let $\sleq E M$ be a sequence of
intervals. Then
$$\sum_{j=1}^M\nm{E_jx}\le C(1+2M/N)\nm x\ .$$

\Proof  For convenience, let us normalize so that $\nm x=N$ and $x=\sm 1 N
x_i$, where $\sleq x N$ and $\nm{x_i}\le C$\ for each $i$.
Given $j$, let $A_j$ be the set of $i$ such that $\supp(x_i)\subset
E_j$ and let $B_j$ be the set of $i$ such that $E_j(x_i)\ne 0$. By the triangle
inequality and the fact that the basis is bimonotone,
$$\nm{E_jx}\le\bgnm{\sum_{i\in B_j}x_i}\le C(|A_j|+2)\ .$$
Since $\sum_{j=1}^M|A_j|\le N$, we get
$$\sum_{j=1}^M\nm{E_jx}\le C(N+2M)$$
which gives the result, because of our normalization.
\hfill $\square$ \bigskip

In order to state the next lemma, we shall need some more definitions. The
first is a technicality. If $f\in\cal F$, let $M_f:\R\ra\R$ be defined by
$M_f(x)=f^{-1}(36x^2)$.

The next definition is of great importance in this paper. We shall say that a
sequence $x_1<\dots<x_N$ is a {\sl rapidly increasing sequence of}
$\ell_{1+}$-{\sl averages}, or R.I.S., {\sl for} $f$ {\sl of length} $N$ {\sl
with constant} $1+\e$ if $x_k$ is an $\ell_{1+}^{n_k}$-average with constant
$1+\e$ for each $k$, $n_1\ge 2(1+\e)M_f(N/\e')/\e' f'(1)$ and
$${\e'\over 2}f(n_k)^{1/2}\ge|\supp(x_{k-1})|$$
for $k=2,\dots,N$. Here $f'(1)$ is the right derivative of $f$ at $1$ and $\e'$
is a useful notation for $\min\{\e,1\}$ which we shall use throughout the
section. Obviously there is nothing magic about the exact conditions in this
definition. The important point is that the $n_k$s increase fast enough, the
speed depending on the sizes of the supports of the earlier $x_j$s.

We make one further definition. A functional $x^*$ is an $(M,g)$-{\sl form} if
$\nm{x^*}\le 1$ and $x^*=\sum_{j=1}^Mx_j^*$ for some sequence
$x_1^*<\dots<x_M^*$ of successive functionals such that $\nm{x_j^*}\le
g(M)^{-1}$ for each $j$.

\proclaim Lemma 3. Let $f,g\in\cal F$, let $g\ge f^{1/2}$ and let $X\in\cal X$
satisfy a lower $f$-estimate. Let $\e>0$, let $\seq x N$ be a R.I.S. in $X$
for $f$ with constant $1+\e$ and let $x=\sum_{i=1}^Nx_i$. Let $M\ge M_f(N/\e')$,
let $x^*$ be an $(M,g)$-form and let $E$ be any interval. Then $|x^*(Ex)|\le
1+\e+\e'$.

\Proof  If $x^*$ is an $(M,g)$-form then so is $Ex^*$ for any interval $E$.
Since $x^*(Ex)=(Ex^*)(x)$, we can forget about the interval $E$ in the
statement of the lemma. For each $i$, let $n_i$ be maximal such that $x_i$ is
an $\ell_{1+}^{n_i}$-average with constant $1+\e$. Let us also write
$x^*=\sum_{j=1}^Mx_j^*$ in the obvious way and set $E_j=\ran(x_j^*)$. We first
obtain three easy estimates for $|x^*(x_i)|$. Since $\nm{x^*}\le 1$, we
obviously have $|x^*(x_i)|\le 1$. Then, since $\nm{x_j^*}\le g(M)^{-1}\le
f(M)^{-1/2}$, we have $|x^*(x_i)|\le f(M)^{-1/2}\sum_{j=1}^M\nm{E_jx_i}$. By
our assumption about $X$, this is at most
$f(M)^{-1/2}f\bigl(|\supp(x_i)|\bigr)$ and by Lemma 2 it is at most
$(1+\e)(1+2Mn_i^{-1})f(M)^{-1/2}$.

Let $t$ be maximal such that $n_t\le M$. Then, if $i<t$, we have
$f\bigl(|\supp(x_i)|\bigr)\le 2^{i-t+1}f\bigl(|\supp(x_{t-1})|\bigr)$, and also
$f\bigl(|\supp(x_{t-1})|\bigr)\le(\e'/2)f(n_t)^{1/2}\le(\e'/2)f(M)^{1/2}$. Using
this and the other two estimates above, we obtain
$$\eqalign{|x^*(x)|&\le\sum_{i=1}^N|x^*(x_i)|\le\e'+1+3(1+\e)(N-t)f(M)^{-1/2}\cr
                   &\le 1+\e'+3(1+\e)N(\e'/6N)\cr
                   &=1+\e'+(\e'/2)(1+\e)\le 1+\e+\e' \cr}$$
as stated. \hfill $\square$ \bigskip

\proclaim Corollary 4. Let $f,X,\e,M,\seq x N$ and $x$ be as in Lemma 3,
let $\sleq E M$ and let $E$ be any interval. Then
$$f(M)^{-1}\sum_{i=1}^M\nm{E_iEx}\le 1+\e+\e'\ .$$

\Proof  Let $x_j^*$ be a support functional of $E_jx$ and let
$x^*=f(M)^{-1}\sum_{i=1}^Mx_i^*$. Then $\nm{x^*}\le 1$ because $X$ satisfies a
lower $f$-estimate. It follows that $x^*$ is an $(M,f)$-form, so we can apply
Lemma 3 with $g=f$. \hfill $\square$ \bigskip

We now introduce a further convenient definition. Let $\sleq x N$ be a R.I.S.
for $f$ with constant $1+\e$, for some $f\in\cal F$ and some $\e>0$. For each
$i$, let $n_i$ be maximal such that $x_i$
is an $\ell_{1+}^{n_i}$-average with
constant $1+\e$ and let us write it out as $x_i=x_{i1}+\dots+x_{i\,n_i}$, where
$\nm{x_{ij}}\le(1+\e)n_i^{-1}$ for each $j$. Given an interval $E\subset\N$,
let $i=i_E$ and $j=j_E$ be respectively minimal and maximal such that $Ex_i$
and $Ex_j$ are non-zero, and let $r=r_E$ and $s=s_E$ be respectively minimal
and maximal such that $Ex_{ir}$ and $Ex_{js}$ are non-zero. Define the {\sl
length} $\l(E)$ of the interval $E$ to be
$j_E-i_E+(s_E/n_{j_E})-(r_E/n_{i_E})$. Thus the length of $E$ is the number of
$x_i$s contained in $E$, allowing for fractional parts. It is easy to check
that if $\sleq E M$ and $E=\bigcup E_i$ then $\sum\l(E_i)\le\l(E)$. Obviously
this definition depends completely on the R.I.S. but it will always be clear
from the context which R.I.S. is being considered.

The next lemma is the most important one for our purposes.

\proclaim Lemma 5. Let $f,g\in\cal F$ with $g\ge\sqrt f$, let $X\in\cal X$
satisfy a lower $f$-estimate, let $\e>0$, let $\sleq x N$ be a R.I.S. in X
for $f$ with constant $1+\e$ and let $x=\sum_{i=1}^Nx_i$. Suppose that
$$\nm{Ex}\le\sup\Bigl\{|x^*(Ex)|:M\ge 2, x^*\ \hbox{is an $(M,g)$-form}
\Bigr\}$$
for every interval $E$ of length at least 1. Then
$\nm x\le(1+\e+\e')Ng(N)^{-1}$.

\Proof  It follows from the triangle inequality that
$\nm{Ex}\le(1+\e)(\l(E)+n_1^{-1})$. If $\l(E)\ge(1+\e)/\e'n_1$ then we get
$\nm{Ex}\le (1+\e+\e')\l(E)$. Let $G$ be defined by $G(x)=x$ when $0\le x\le 1$
and $G(x)=xg(x)^{-1}$ when $x\ge 1$. Recall that, because of the properties of
the set $\cal F$, $G$ is concave and increasing on $[1,\infty)$ and satisfies
$G(xy)\ge G(x)G(y)$ for every $x,y$ in the same interval. It is easy to check
that these properties are still true on the whole of $\R_+$. We shall show that
if $\l(E)\ge (1+\e)/\e'n_1$, then $\nm{Ex}\le(1+\e+\e')G(\l(E))$. The remarks
we have just made show this when $\l(E)\le 1$.

Let us suppose then that $E$ is a minimal interval of length at least
$(1+\e)/\e' n_1$ for which the inequality fails. We know that $\l(E)\ge 1$. We
also know that there exists some $(M,g)$-form $x^*=\sum_{i=1}^Mx_i^*$ such that
$\nm{Ex}\le|x^*(Ex)|$. By Lemma 3, we must have $M\le M_f(N/\e')$ or the
inequality would not fail for $E$. Letting $E_i=E\cap\ran(x_i^*)$, we have
$$\nm{Ex}\le g(M)^{-1}\sum_{i=1}^M\nm{E_ix}$$
by the definition of an $(M,g)$-form.

Let $\l_i=\l(E_i)$ for each $i$. For each $i$ we either have $\l_i\le(1+\e)/\e'
n_1$ or, by the minimality of $E$, that $\nm{E_ix}\le(1+\e+\e')G(\l_i)$. Let
$A$ be the set of $i$ with the first property and let $B$ be the complement of
$A$. Let $k$ be the cardinality of $A$.

Since $G$ is a concave and non-decreasing function and $\sum\l_i\le\l$,
Jensen's inequality gives us that
$$\eqalign{\sum_{i\in B}\nm{E_ix}&\le(1+\e+\e')\sum_{i\in B}G(\l_i)\cr
                                 &\le(1+\e+\e')(M-k)G(\l/(M-k))\ .\cr}$$
It follows that
$$\eqalign{\nm{Ex}
       &\le M^{-1}G(M)\Bigl[(1+\e+\e')(M-k)G(\l/(M-k))+
	   (1+\e)(1+\e+\e')k/\e'n_1\Bigr]\cr
       &\le (1+\e+\e')\Bigl[(1-k/M)G(M)G(\l/(M-k))+(1+\e)k/\e'n_1\Bigr]\cr
 &\le (1+\e+\e')\Bigl[(1-k/M)G\bigl((1-k/M)^{-1}\l\bigr)+(1+\e)k/\e'n_1\Bigr]
                                                                       \ .\cr}$$

Let $G'(1)$ be the right derivative of $G$ at $1$. Since $G$ is a concave
function we have the easy inequality
$$(1-t)G\bigl({\l\over 1-t}\bigr)+t(G(1)-G'(1))\le G(\l)$$
for every $0\le t<1$ and $\l\ge 1$. Also, $G(1)-G'(1)=1-G'(1)=g'(1)$, and since
$g\ge\sqrt f$ we have $g'(1)\ge f'(1)/2>0$.

By the definition of R.I.S. we have $n_1\ge 2(1+\e)M_f(N/\e')/\e'f'(1)$. It
follows that
$$(1+\e)k/\e'n_1\le(k/M)\bigl((1+\e)M_f(N/\e')/\e'n_1\bigr)
\le(k/M)(f'(1)/2)\le(k/M)g'(1)\ .$$
Hence, by the inequality above with $t=k/M$, we have
$$\eqalign{(1+\e+\e')\Bigl[(1-k/M)G\bigl((1-k/M)^{-1}\l\bigr)&+(1+\e)k/\e'n_1
\Bigr] \cr
&\le (1+\e+\e')\Bigl((1-t)G\bigl({\l\over 1-t}\bigr)+tg'(1)\Bigr) \cr
&\le (1+\e+\e')G(\l)\ ,\cr}$$
contradicting our assumption about the interval $E$, and proving the lemma.
\hfill $\square$ \bigskip

It is now easy to construct an asymptotic biorthogonal system in Schlumprecht's
space. Let $\delta \in (0,1)$, let $N_1<N_2<\dots$ be
a sequence of integers satisfying
$f(N_1)/N_1<\d/2$, $f(N_1)>8\d^{-1}$ and $N_j>M_f(2N_{j-1})$ for all $j>1$.
Let $A_k$ be the set of norm-1 vectors of the form $x=\sum_{i=1}^{N_k}x_i$
where $x_1,\dots,x_{N_k}$ is a multiple of a R.I.S. with constant $1+\d/2$.
Because Schlumprecht's space satisfies a lower $f$-estimate, we know that the
multiple is at most $f(N_k)N_k^{-1}$. Let $A_k^*$ be the set of functionals of
the form $f(N_k)^{-1}\sum_{i=1}^{N_k}x_i^*$ where $x_1^*<\dots<x_{N_k}^*$ and
$\nm{x_i^*}\le 1$ for each $i$. It is clear that the sets $A_k$ are asymptotic
for every $k$. If $j>k$ then using the fact that $N_j>M_f(2N_k)$, we may apply
Lemma 3 with $\e=1/2$ and $M=N_j$ since $y^*$ is clearly an $(M,f)$-form
whenever $y^*\in A_j^*$. Because of the normalization of the R.I.S., this gives
us $|y^*(x)|\le 2f(N_k)/N_k < \delta$
for every $y^*\in A_j^*$ and $x\in A_k$.

If $j<k$ then we know from Lemma 5 that $\nm{\sum_{i\in A}x_i}\le
2|A|f(N_k)/N_kf(|A|)$ for every subset $A$ of $\{1,2,\dots,N_k\}$. If
$|A|\ge\sqrt{N_k}$ then this is at most $4|A|/N_k$. By splitting into
$\sqrt{N_k}$ successive pieces of this form, we find that $x$ is an
$\ell_{1+}^{\sqrt{N_k}}$-average with constant~4. By Lemma 2 we obtain that
$|y^*(x)|\le f(N_j)^{-1}.4\bigl(1+2N_j/\sqrt{N_k}\bigr)\le
8f(N_j)^{-1} < \delta$.

Finally, we know that $\nm x\le (1+\d)N_kf(N_k)^{-1}\nm{x_i}$ for each $i$, so
if we let $x_i^*$ be a support functional of $x_i$ then
$x^*=f(N_k)^{-1}\sum_{i=1}^{N_k}x_i^*$ is an element of $A_k^*$ and
$x^*(x)\ge(1+\d)^{-1} > 1-\d$. It follows that $A_1,A_2,\dots$ and
$A_1^*,A_2^*,\dots$ form an asymptotic biorthogonal system with constant $\d$.

This together with the result of the last section shows that, for every $C$,
Schlumprecht's space can be renormed so as not to contain a $C$-unconditional
basic sequence. Since Schlumprecht's space itself has a 1-unconditional basis,
it follows that it is arbitrarily distortable. This is also an easy
direct consequence of the existence of an asymptotic biorthogonal system, or
indeed from Lemma 5, which is what Schlumprecht used.

\bigbreak

\noindent {\bf \S 3. A space containing no unconditional basic sequence.}

We now come to the main result of the paper, namely the construction of a
Banach space $X$ containing no unconditional basic sequence. As we mentioned in
the introduction, it was observed by W. B. Johnson that our original arguments
could be modified to show that $X$ was actually H.I. This is what we shall
actually present in this section.

The definition of the space resembles that of Schlumprecht's space, or at least
can do. We shall actually give three equivalent definitions, for which we shall
need a certain amount of preliminary notation.

First, let ${\bf Q}$ be the set of real sequences with finite support, rational
coordinates and maximum at most 1 in modulus. Let $J\subset\N$ be a set such
that, if $m<n$ and $m,n\in J$, then $\log\log\log n\ge 2m$. Let us write $J$ in
increasing order as $\{j_1,j_2,\dots\}$. We shall also assume that $f(j_1)\ge
36$. (Recall that $f(x)$ is the function $\log_2(x+1)$). Now let $K\subset J$ be
the set $\{j_2,j_4,j_6,\dots\}$ and let $L\subset\N$ be the set of integers
$j_1,j_3,j_5,\dots$.

Let $\sigma$ be an injection from the collection of finite sequences of
successive elements of ${\bf Q}$ to $L$ such that, if $\seq z s$ is such a
sequence, $S=\sigma(\seq z s)$ and $z=\sm i s z_i$ then $(1/20)f\Bigl(S^{1/40}
\Bigr)^{1/2}\ge |\supp(z)|$.

We shall use the injection $\sigma$ to define special
functionals in an arbitrary normed space of the form $(c_{00},\nm.)$ in much
the same way that we defined them in Section 1. (Of course, for most
spaces they are not terribly interesting).

If $X=(c_{00},\nm.)$ is a normed space on the finitely supported sequences, and
$m\in\N$, let $A_m^*(X)$ be the set of functionals of the form $f(m)^{-1}
\sum_{i=1}^mf_i$ such that $\sleq f m$ and $\nm{f_i}\le 1$ for each $i$. If
$k\in\N$, let $\G_k^X$ be the set of sequences $\seq g k$ such that $g_i\in
{\bf Q}$ for each $i$, $g_1\in A_{j_{2k-1}}^*(X)$ and $g_{i+1}\in
A_{\sigma(\seq g i)}^*(X)$ for each $1\le i\le k-1$. We call these {\sl special
sequences}. Let $B_k^*(X)$ be the set of functionals of the form
$f(k)^{-1/2}\sum_{j=1}^kg_j$ such that $(\seq g k)\in\G_k^X$. These are {\sl
special functionals}.

Dually, if a convex set $D\subset c_{00}$ is given, we define $A_m(D)$ to be
the set of vectors of the form $f(m)^{-1}\sum_{i=1}^mx_i$ where
$\sleq x m$ and $x_i\in D$ for each $i$. Then special sequences of vectors are
defined using $\sigma$ in the obvious corresponding way, and this gives us sets
$B_k(D)$.

Our first definition of the norm is geometrical, and goes via the dual space.
Let $D_0$ be the intersection of $c_{00}$ with the unit ball of $\ell_1$. Once
we have defined $D_N$, let $D_N'$ be the set of vectors of the form
$f(N)^{-1}\sum_{i=1}^Nx_i$, where $\seq x N$ are successive vectors in $D_N$.
Let $D_N''$ be the set of special vectors for $D_N$ with lengths in $K$, that
is, $D_N''=\bigcup_{k\in K}B_k(D_N)$. Let $D_N'''$ be the set of vectors $Ex$
where $x\in D_N$. Then let $D_{N+1}$ be the convex hull of the union of $D_N'$,
$D_N''$ and $D_N'''$.

Now let $D=\bigcup_{N=0}^\infty D_N$. It is easy to see that $D$ is the
smallest convex set closed under taking sums of the form
$f(N)^{-1}\sum_{i=1}^Nx_i$, taking special vectors with lengths in $K$ and
under interval projections. Our space is defined by
$$\nm x=\sup\{|\sp{x,y}|:y\in D\}\ .$$

The second definition of the norm is as the limit of a sequence of norms.
Define $X_0=(c_{00},\nm._0)$ by $\nm x_0=\nm x_\infty$, and, for $n\ge 0$, let
$$\eqalign{\nm x_{X_{n+1}}&=
\sup\Bigl\{f(n)^{-1}\sum_{i=1}^n\nm{E_ix}_{X_n}:n\in\N,\sleq E n\Bigr\}\cr
&\vee\sup\Bigl\{|g(Ex)|:k\in K,\, g\in B_k^*(X_n),\, E\subset\N\Bigr\}\ .
\cr}$$
Note that this is an increasing sequence of norms, because the sets
$B_k^*(X_n)$ increase as $n$ increases (and more generally, if $\nm x_Y\le\nm
x_Z$ for every $x\in c_{00}$, then $B_k^*(Y)\subset B_k^*(Z)$). They are also
all bounded above by the $\ell_1$-norm. Define $\nm.$ by $\nm
x=\lim_{n\ra\infty}\nm x_{X_n}$.

Finally, we give an implicit definition of the norm in the obvious way. Set
$$\eqalign{\nm x=\nm x_{c_0}&\vee
\sup\Bigl\{f(n)^{-1}\sum_{i=1}^n\nm{E_ix}:2\le n\in\N,\sleq E n\Bigr\}\cr
&\vee\sup\Bigl\{|g(Ex)|:k\in K,\, g\in B_k^*(X),\, E\subset\N\Bigr\}\ .
\cr}$$

Recall that $E\subset\N$ is always an interval in these definitions. Its role
is to ensure that $(\bfe_i)_{i=1}^\infty$ is a (bimonotone) normalized Schauder
basis for the completion of $X$. Note also that if we did not insist that the
$E_i$ were intervals then the unit vector basis of this space would trivially
be unconditional. It is not hard to check that the norm given by the third
definition is indeed well-defined and agrees with both the previous ones.

Before getting down to analysing the space, we shall need a few simple facts
about functions in the class $\cal F$ defined earlier.

We shall now introduce some convenient definitions. Let
$f:[1,\infty)\ra[1,\infty)$ be a function. The (increasing) {\sl
submultiplicative hull} of $f$ is the function $F$ defined by
$$F(x)=\inf\Bigl\{f(x_1)f(x_2)\dots f(x_k):k\in\N, x_i\ge 1, x_1\dots x_k\ge x
\Bigr\}\ .$$
The following facts are trivial. First, $F\le f$. Second, $F(xy)\le F(x)F(y)$.
Third, if $g:[1,\infty)\ra[1,\infty)$ is any non-decreasing submultiplicative
function dominated by $f$, then $g$ is dominated by $F$. (That is, $F$ is the
largest non-decreasing submultiplicative function dominated by $f$).

Now let $g:[1,\infty)\ra[1,\infty)$ be any function. The {\sl concave envelope}
of $g$ is of course the smallest concave function $G:[1,\infty)\ra[1,\infty)$
dominating $g$, that is,
$$G(x)=\sup\Bigl\{\l g(y)+(1-\l)g(z):0\le\l\le 1, \l y+(1-\l)z=x\Bigr\}\ .$$

We now prove an easy lemma.

\proclaim Lemma 6. If $g:[1,\infty)\ra[1,\infty)$ is a supermultiplicative
function, then its concave envelope is also supermultiplicative.

\Proof   Let $\e>0$ and let $x_1,x_2\ge 1$. We shall show that
$(G(x_1)-\e)(G(x_2)-\e)\le G(x_1)G(x_2)$, which will prove the result.
First, for $i=1,2$, pick $\l_i,\mu_i,y_i$ and $z_i$ such that $0\le\l_i\le 1$,
$\l_i+\mu_i=1$, $\l_iy_i+\mu_iz_i=x_i$ and $\l_ig(y_i)+\mu_ig(z_i)\ge
G(x_i)-\e$.

Then
$$\eqalign{(G(x_1)-\e)(G(x_2)-\e)
&\le(\l_1g(y_1)+\mu_1g(z_1))(\l_2g(y_2)+\mu_2g(z_2)) \cr
&\le\l_1\l_2g(y_1y_2)+\l_1\mu_2g(y_1z_2)+\mu_1\l_2g(z_1y_2)+\mu_1\mu_2g(z_1z_2)
\cr
&\le\l_1G(y_1x_2)+\mu_1G(z_1x_2)\le G(x_1x_2)\cr}$$
as we wanted. \hfill $\square$ \bigskip

Now let us define a function $\phi:[1,\infty)\ra[1,\infty)$ as follows.
$$\phi(x)=\cases{(\log_2(x+1))^{1/2}&if $x\in K$ \cr
                  \log_2(x+1)        &otherwise. \cr}$$
Let $h$ be the submultiplicative hull of $\phi$, let $H$ be the function given
by $H(x)=x/h(x)$ and let $G$ be the concave envelope of $H$. Since $h$ is
submultiplicative, $H$ is supermultiplicative, so $G$ is also
supermultiplicative. Now let $g(x)=x/G(x)$. Then $g$ is submultiplicative. As
before, let $f$ be the function $\log_2(x+1)$. The easy facts about
submultiplicative hulls and concave envelopes and the fact that
$\sqrt f \in {\cal F}$ give that $(\log_2(x+1))^{1/2}\le
g(x)\le\phi(x)\le\log_2(x+1)$. It will also be useful to extend the definition
of $G$ to the whole of $\R_+$ by setting $G(x)=x$ when $0\le x\le 1$. It is easy
to check that $G$ thus extended is still supermultiplicative and concave, as we
commented in the proof of Lemma 5.

We now need to calculate $G(N)$ when $N\in L$. In fact we shall want slightly
more than this, as is suggested by the statement of the next lemma.

\proclaim Lemma 7. If $N\in L$ then $G(x)=xf(x)^{-1}$ for every $x$ in the
interval $[\log N,\exp N]$.

\Proof  Let $k,l\in K$ be maximal and minimal respectively such that $k<N$ and
$l>N$, and let $(k!)^4<x<f^{-1}(f(l)^{1/2})$. We shall show first that
$h(x)=f(x)$.
$h(x)$ is defined to be $\inf\Bigl\{\phi(x_1)\dots\phi(x_m):x_i\ge
1,x_1x_2\dots x_m\ge x\Bigr\}$. We know that $\phi(x_j)=f(x_j)^{1/2}$ if
$x_j\in K$ and $f(x_j)$ otherwise. By the submultiplicativity of $f$, there can
be at most one $j$ such that $x_j>1$ and $x_j\notin K$. Because
$f(l)^{1/2}>f(x)$ we also know that if $x_j\in K$ then $x_j\le k$. Thirdly, it
is not possible to find $r,s,t$ such that $x_r=x_s=x_t\in K$ since, for every
$p\in K$, $f(p)^{3/2}>f(p^3)$ by our choice for $\min(K)$. Since $x>(k!)^4$ it
is clear that at least one, and hence exactly one $x_j$ is not in $K$. Let it
be $x_1$ and assume that $m>1$. Now we know that $x_2x_3\dots x_m\le (k!)^2\le
x^{1/2}$. It follows that $x_1\ge x^{1/2}$ and hence that
$\phi(x_1)\dots\phi(x_m)\ge f(x^{1/2})f(\min(K))^{1/2}$. Since $f$ is the
function $\log_2(x+1)$ and $f(\min(K))\ge 36$, this is greater than $f(x)$.
This contradiction shows that $m=1$, hence $h(x)=f(x)$.

This shows that $H(x)=xf(x)^{-1}$ whenever $(k!)^4<x<f^{-1}(f(l)^{-1/2})$, and
in particular for all $x$ in the interval $[\log\log N,\exp\exp N]$. It is easy
to deduce from this the conclusion of the lemma. Indeed, given $x_0$ in the
interval $[\log N,\exp N]$, we will certainly know that $G(x_0)=x_0f(x_0)^{-1}$
if the function given by the tangent to $xf(x)^{-1}$ at $x_0$ is at least
$xf(x)^{-1/2}$ for all positive $x$ outside the interval $[\log\log N,\exp\exp
N]$.

The equation of the tangent at $x_0$ is
$$y={x_0\over f(x_0)}+{1\over f(x_0)}\left(1-{x_0\over(x_0+1)\log(x_0+1)}\right)
(x-x_0)\ .$$

When $x\ge 0$ this is certainly at least $x_0^2\log 2/(x_0+1)(\log(x_0+1))^2$
which is at least $x_0/2f(x_0)^2$. For $x_0\ge\log N$ and $x\le\log\log N$ this
exceeds $xf(x)^{-1/2}$. When $x\ge 2x_0$ we also know that $y\ge x/4f(x_0)$.
When $x\ge\exp\exp N$ the condition $x_0\le\exp N$ is enough to guarantee that
this is at least as big as $xf(x)^{-1/2}$.
\hfill$\square$\bigskip

We shall now prove a crucial lemma about $X$. It is an easy consequence of
Lemmas~5 and~7.

\proclaim Lemma 8. Let $N\in L$, let $n\in[\log N,\exp N]$, let $\e>0$ and
let $\seq x n$ be a R.I.S. with constant $1+\e$. Then
$\nm{\sum_{i=1}^nx_i}\le(1+\e+\e')nf(n)^{-1}$.

\Proof  It is obvious from the implicit definition of the norm in $X$ that it
satisfies a lower $f$-estimate. Let $g$ be the function defined before the last
lemma. As usual let $x=\sum_{i=1}^nx_i$.
Since $g\le\phi$, it is clear that every vector in $X$ is either normed by
an $(M,g)$-form or has the supremum norm. It is also clear that the second
possibility does not happen in the case of vectors of the form $Ex$ when
$\l(E)\ge 1$. Since $g\in\cal F$ and, as we commented above, $g\ge f^{1/2}$,
all the hypotheses of Lemma 5 are satisfied. It follows that
$\nm{\sum_{i=1}^nx_i}\le(1+\e+\e')G(n)$. By Lemma 7, $G(n)=nf(n)^{-1}$, so the
lemma is proved. \hfill$\square$\bigskip

\proclaim Lemma 9. Let $N\in L$, let $0<\e<1/4$, let $M=N^{\e}$
and let $\seq x N$ be a R.I.S. with constant $1+\e$. Then $\sm i Nx_i$ is an
$\ell_{1+}^M$-vector with constant $(1+4\e)$.

\Proof  Let $m=N/M$, let $x=\sum_{i=1}^Nx_i$ and for $1\le j\le M$ let
$y_j=\sum_{i=(j-1)m+1}^{jm}x_i$. Then each $y_j$ is the sum of a R.I.S. of
length $m$ with constant $(1+\e)$. By Lemma 8 we have
$\nm{y_j}\le(1+2\e)mf(m)^{-1}$ for every $j$ while $\nm{\sum_{j=1}^my_j}=\nm
x\ge Nf(N)^{-1}$. It follows that $x$
is an $\ell_{1+}^M$-vector with constant at
most $(1+2\e)f(N)/f(m)$. But $m=N^{1-\e}$ so $f(N)/f(m)\le(1-\e)^{-1}$. The
result follows. \hfill $\square$ \bigskip

We shall now prove that $X$ is H.I. As we noted earlier, this implies that $X$
contains no unconditional basic sequence, but in proving that $X$ is H.I. we
shall more or less have proved that directly anyway.

Let $Y,Z$ be two infinite-dimensional subspaces of $X$ such that $Y\cap
Z=\{0\}$. Our aim is now to show that the projection from $Y+Z$ to $Y$ given by
$y+z\mapsto y$ is not continuous. To do this, we shall construct, for every
$\d>0$, vectors $y\in Y$ and $z\in Z$ of norm at least~1 such that
$\nm{y-z}<\d$. This implies that the above projection has norm at least
$\d^{-1}$, proving the result. So let us now choose $\d>0$ and let $k\in K$ be
an integer such that $f(k)^{-1/2}<\d$.

By standard arguments, we may assume that both $Y$ and $Z$ are spanned by block
bases. Since $X$ satisfies a lower $f$-estimate, Lemma 1 tells us that every
block subspace of $X$ contains, for every $\e>0$ and $N\in\N$, an
$\ell_{1+}^N$-average with constant $1+\e$. It is also immediate from the
definition of the norm that every vector either has the supremum norm or
satisfies the inequality
$$\nm{Ex}\le\sup\Bigl\{|x^*(Ex)|:M\ge 2, x^*\ \hbox{is an $(M,g)$-form}
\Bigr\}$$
where $g$ is the function obtained from $\phi$ after the proof of Lemma 6. This
allows us to make the following construction.

Let $x_1\in Y$ be a normalized R.I.S. vector of length $M_1=j_{2k-1}\in L$ and
constant $(1+\e/4)$,
where $\e=1/10$ and $M_1^{\e/4}\ge N_1=4M_f(k/\e)/\e f'(1)$.

Let the normalized R.I.S. whose sum is $x_1$ be $x_{11},\dots,x_{1M_1}$. By
Lemma 9, $x_1$ is an $\ell_{1+}^{N_1}$-average
with constant $1+\e$. By Lemma 5,
we have $\nm{x_1}\le(1+\e)M_1g(M_1)^{-1}\nm{x_{11}}$. For each $j$ between 1
and $M_1$ let $x_{1j}^*$ be a support functional for $x_{1j}$ and let
${x'}_1^*$ be the $(M_1,g)$-form $g(M_1)^{-1}\sum_{j=1}^{M_1}x_{1j}$. Then
${x'}_1^*(x_1)\ge(1+\e)^{-1}\nm{x_1}$. By continuity and the density of ${\bf
Q}$ it follows that there exists an $(M_1,g)$-form $x_1^*\in {\bf Q}$ such that
$|x_1^*(x_1)-1/2|\le k^{-1}$ and $\ran(x_1^*)=\ran(x_1)$. Also, note that by
Lemma 7 there is no difference between an $(M_1,g)$-form and an $(M_1,f)$-form.

Now let $M_2=\sigma(x_1^*)$ and pick a normalized R.I.S. vector $x_2\in Z$ of
length $M_2$ with constant $1+\e/4$ such that $x_1<x_2$. Then $x_2$ is an
$\ell_{1+}^{N_2}$-average with constant $1+\e$,
where $N_2=M_2^{\e/4}$. As above,
we can find an $(M_2,g)$-form $x_2^*$ such that $|x_2^*(x_2)-1/2|\le k^{-1}$
and $\ran(x_2^*)=\ran(x_2)$.

Continuing in this manner, we obtain a pair of sequences $\seq x k$ and
$x_1^*,\dots,x_k^*$ with various properties we shall need.
First, $x_i\in Y$ when $i$ is odd and $Z$ when $i$ is even.
Second, $\nm{x_i}=1$
for every $i$ and $\nm{x_i^*}\le 1$. We also know that $|x_i^*(x_i)-1/2|\le1/k$
for each $i$. Recall that $\sigma$ was chosen so that if $\seq z s$ is a
sequence of successive vectors in ${\bf Q}$, $S=\sigma(\seq z s)$ and $z=\sm i
s z_i$ then $(1/20)f\Bigl(S^{1/40}\Bigr)^{1/2}\ge |\supp(z)|$.

This and the lower
bound for $N_1$ ensure that $\seq x k$ is a R.I.S. of length $k$ with constant
$1+\e$. It will also be important that the sequence
$x_1,-x_2,x_3,\dots,(-1)^{k-1}x_k$ is a R.I.S. which is obviously true as well.
Finally, and perhaps most importantly, the sequence $x_1^*,\dots,x_k^*$ has
been carefully chosen to be a special sequence of length $k$. It follows
immediately from the implicit definition of the norm and the fact that
$\ran(x_i^*)\subset\ran(x_i)$ for each $i$ that
$$\nm{\sum_{i=1}^kx_i}\ge f(k)^{-1/2}\sum_{i=1}^kx_i^*(x_i)\ge
f(k)^{-1/2}(k/2-1)\ .$$

The proof will be complete if we can find a suitable upper bound for
$\nm{\sum_{i=1}^k(-1)^{i-1}x_i}$. To do this we shall apply Lemma 5 one final
time. First, we shall show that if $z_1^*,\dots,z_k^*$ is any special sequence
of functionals of length $k$ and $E$ is an interval of length at least~1 with
respect to the R.I.S. $x_1,-x_2,\dots,(-1)^{k-1}x_k$, then $|z^*(Ex)|\le 1/2$,
where $z^*$ is the $(M,g)$-form $f(k)^{-1/2}\sum_{i=1}^kz_i^*$ and
$x=\sum_{i=1}^k(-1)^{i-1}x_i$.

Indeed, let $t$ be maximal such that $z_t^*=x_t^*$ or zero if no such $t$
exists. We know that if $i\le t$ then $|z_i^*(x_i)-1/2|\le k^{-1}$.
Suppose $i\ne j$ or one of $i,j$ is greater than $t$. Then since $\sigma$ is an
injection, we can find $L_1\ne L_2\in L$ such that $z_i^*$ is an $(L_1,f)$-form
and $x_j$ is the normalized sum of a R.I.S. of length $L_2$ and also an
$\ell_{1+}^{L'_2}$-average with constant $1+\e$, where $L'_2=L_2^{\e/4}$. Just
as at the end of section 2, we can now use Lemmas 2 and 5 to show that
$|z_i^*(x_j)|<k^{-2}$.

It follows that
$$\eqalign{\Bigl|\Bigl(\sum_{i=1}^kz_i^*\Bigr)(x)\Bigr|&\le
k^2.k^{-2}+\Bigl|\sum_{i=1}^t(-1)^{i+1}z_i^*(x_i)\Bigr|
+|z_{t+1}^*(x_{t+1})|+\sum_{i=t+2}^k|z_i^*(x_i)|  \cr
&\le 1+(1+k.k^{-1})+1+k^2.k^{-2}\le 5\ . \cr}$$

The interval $E$ is easily seen to increase this to at most 6. It follows that
$|z^*(Ex)|\le 6f(k)^{-1/2}<1/2$ as claimed.

Now let $\phi'$ be the function
$$\phi'(x)=\cases{(\log_2(x+1))^{1/2}&if $x\in K, x\ne k$ \cr
                  \log_2(x+1)       &otherwise.          \cr}$$
Let $g'$ be the function obtained from $\phi'$ just as $g$ was obtained from
$\phi$. It follows from our remarks about special sequences of length $k$ that
$$\nm{Ex}\le\sup\Bigl\{|x^*(Ex)|:M\ge 2, x^*\ \hbox{is an $(M,g')$-form}
\Bigr\}$$
whenever $E$ is an interval of length at least 1. Since $x$ is the sum of a
R.I.S. it follows that we can use Lemma 5 to show that
$\nm{x}\le(1+2\e)kg'(k)^{-1}$. Finally, Lemma 7 gives us that $g'(k)=f(k)$,
since $\phi\le\phi'$ yields $g\le g'\le f$.

We have now constructed two vectors $y\in Y$, the sum of the odd-numbered
$x_i$s, and $z\in Z$, the sum of the even-numbered $x_i$s, such that
$\nm{y+z}\ge(1/3)f(k)^{1/2}\nm{y-z}$. Hence $Y$ and $Z$ do not form a
topological direct sum, so $X$ is H.I. If $X$ contained an unconditional basic
sequence $x_1,x_2,\dots$ then the subspace generated by this sequence would
split into a direct sum of the subspaces generated by $\{x_{2n-1}:n\in\N\}$ and
$\{x_{2n}:n\in\N\}$. It follows that $X$ does not contain an unconditional
basic sequence. The reader will observe that it is easy to use the preceding
argument to show this directly. In the next section, we shall examine some of
the other consequences of a space being H.I., but first we shall observe that
$X$ is reflexive. For the definitions of the terms ``shrinking'' and
``boundedly complete'' see [LT, section 1.b].

First, it follows immediately from the fact that $X$ satisfies a lower
$f$-estimate that the standard basis $\bfe_1,\bfe_2,\dots$ is boundedly
complete. Now suppose that it is not a shrinking basis. Then we can find
$\e>0$, a norm-1 functional $x^*\in X^*$ and a sequence of
normalized blocks
$x_1,x_2,\dots$ such that $x^*(x_n)\ge\e$ for every $n$. It follows that
$\sum_{x\in A}x_n$ is an $\ell_{1+}^{|A|}$-vector with constant $\e^{-1}$
for every $A\subset\N$. Given $N\in L$ we may construct a R.I.S. $\seq y N$
with constant $\e^{-1}$ where $y_i$ is of the form $\l_i\sum_{j\in
A_i} x_j$, with $\l_i \ge |A_i|^{-1}$.
Then $x^*(\speq y N)\ge \e N$. For $N$ sufficiently large,
this contradicts Lemma 8.

\bigbreak

\noindent{\bf \S 4. Operators on H.I. spaces.}
\medskip

In this section, we shall prove some results about H.I. spaces over $\CC$. This
is because we shall need to use a little spectral theory. In the next section
we shall show that some of the results carry over to the real case. We do not
know of a direct proof.

Let $X$ be a complex Banach space and let $T$ be a bounded linear operator from
$X$ into itself. We say that $\lambda\in \CC$ is {\sl infinitely singular for}
$T$ if, for every $\e>0$, there exists an infinite-dimensional subspace $Y_\e$
of $X$ such that the restriction of $T-\l I$ to $Y_\e$ has norm at most $\e$.

Saying that $\l$ is {\it not} infinitely singular for $T$ is equivalent to
saying that $T-\l I$ is an isomorphism on some finite-codimensional subspace
of $X$. Since this property is clearly unaffected by a small enough
perturbation, it follows that
$$F_T=\{\l\in\CC:\l\hbox{\ not infinitely singular for\ } T\}$$
is an open subset of $\CC$. Notice that $\ker(T-\l I)$ is finite dimensional
when $\l\in F_T$. We shall now prove some lemmas about $F_T$.

\proclaim Lemma A. If $\l\in F_T$ and if $(x_n)$ is a bounded sequence such
that $(T-\l I)x_n$ is norm-convergent, then $(x_n)$ has a norm-convergent
subsequence; furthermore, the image by $T-\l I$ of any closed subspace of
$X$ is closed.

\Proof  Let $S=T-\l I$, let $Y$ be a finite-codimensional subspace on which $S$
is an isomorphism and let $X=Y\oplus Z$. Let $x_n=y_n+z_n$ with $y_n\in Y$ and
$z_n\in Z$. Then $Sx_n=Sy_n+Sz_n$. Since $Z$ is finite-dimensional and $(x_n)$
is bounded, we can pass to a subsequence such that $Sz_n$ converges. Since
$Sx_n$ converges this gives us that $Sy_n$ converges (relabelling the
subsequence as $Sy_n$). Since $S$ is an isomorphism on $Y$ it follows that
$y_n$ converges. Finally pass to a further subsequence on which $z_n$
converges. To prove the second assertion, note that if $F$ is a closed subspace
of $X$, then $F= F\cap Y + G$, for some finite-dimensional $G$, and hence $T(F)
= T(F\cap Y) + T(G)$ is closed. \hfill $\square$ \bigskip

\proclaim Lemma B. If $\l\in\partial Sp(T) \cap F_T$, then $\l$ is an
eigenvalue of $T$ with finite multiplicity.

\Proof  Since $\l\in\partial Sp(T)$ it is an approximate eigenvalue of $T$.
Hence, there exists a sequence $x_n$ of norm-one vectors with $Tx_i-\l
x_i\rightarrow 0$. By the previous lemma it has a convergent subsequence. But
then the limit of the subsequence is an eigenvector with eigenvalue $\l$.
\hfill $\square$ \bigskip

The next lemma follows easily from well known facts in Fredholm
theory. The argument here is elementary. It was shown to
us by W. B. Johnson, as was the proof of Lemma~D.

\proclaim Lemma C. If $\l\in\partial Sp(T) \cap F_T$, then $\l$ is an isolated
point of $Sp(T)$.

\Proof  Since $F_T$ is open it is enough to show that $\l$ is an isolated point
of $\partial Sp(T) \cap F_T$. Suppose that this is not the case. Then there
exists a sequence $(\l_n)$ in $\partial Sp(T)\cap F_T$ converging to $\l$, with
$\l_n\ne\l$ for every $n$. Since  $\l_n\in F_T$,  $\l_n$ is an eigenvalue, by
Lemma~B. Let $x_n$ be a norm one eigenvector with eigenvalue $\l_n$. By
Lemma~A, since $(T-\l I)x_n$ tends to $0$, we may assume that $(x_n)$ is
norm-convergent to some (norm one) vector $x$ such that $Tx=\l x$. Let $Y$ be
the closed subspace of $X$ generated by the sequence $(x_n)$. Let $U$ be the
restriction of $T-\l I$ to $Y$. It is clear that $Y$ is invariant under $U$ and
that $UY$ is dense in $Y$. Furthermore, since $(T-\l I)Y = UY$ and $\l\in F_T$,
it follows from Lemma~A that $UY$ is closed, and hence that $UY=Y$. Since $x\in
Y$, we know that $Y_0=\ker U$ is not $\{0\}$, and that it is
finite-dimensional. We can therefore write $Y$ as a direct sum $Y_0+Y_1$. We
have that $UY_1=Y$, so for small $\e$ it is still true that $(U-\e I)Y_1=Y$.
But since $(U-\e I)Y_0=Y_0$ when $\e\ne 0$, this yields that $\ker(U-\e
I)\ne\{0\}$, for every small $\e$, contradicting the fact that $\l\in\partial
Sp(T)$. \hfill $\square$ \bigskip

\proclaim Lemma D. Let $Y$ be a subspace invariant under $T$, let $S$ be the
restriction of $T$ to $Y$, and suppose that $Sp(S)=\{\l\}$. If $\l$ is not
infinitely singular for $S$, then $Y$ is finite-dimensional.

\Proof  Suppose that $\l\in F_S$ but that $Y$ is infinite-dimensional. Then
$U=S-\l I$ is an isomorphism on some finite-codimensional subspace $Z$ of $Y$,
and $Sp(U)=\{0\}$. Replacing $U$ by an appropriate multiple, we may assume that
$\|Uz\|\ge\|z\|$ for every $z\in Z$. Define $Z_0=Z$, $Z_1=Z\cap UZ,
\ldots,Z_{k+1}=Z\cap UZ_k$. All these subspaces of $Y$ are
infinite-dimensional. If $z$ is a non-zero element of $Z_k$, we see that $z=
U^k z_0$ for some $z_0\in Z$ and $0<\nm{z_0}\le\nm{U^kz_0}$. This shows that
$\nm{U^k}\ge 1$ for every $k$, contradicting the fact that the spectral radius
of $U$ is $0$. \hfill $\square$ \bigskip

Suppose now that $X$ is a complex H.I. Banach space. Let $T$ be a bounded
linear operator from $X$ into itself. It follows easily from the H.I.
property that there exists at most one value $\l_0$ that is infinitely
singular for $T$. If $\l_0$ is infinitely singular for $T$, the H.I. property
implies that $T-\l_0 Id$ is not an isomorphism on any infinite-dimensional
subspace of $X$.  In other words, $T-\l_0 I$ is strictly singular.

It follows from Lemma~C that the spectrum of $T$ is finite, or consists of a
sequence of eigenvalues converging to $\l_0$. In the second case, it is clear
that $\l_0$ is infinitely singular for $T$. We must check that some $\l$ is
infinitely singular for $T$ in the case of a finite spectrum.

Assume then that $\l\in \partial Sp(T)\cap F_T$. Then $\l$ is isolated in
$Sp(T)$, by Lemma~C. If $Q$ is the spectral projection associated with $\l$,
then $Y = QX$ is finite-dimensional, by the spectral mapping theorem and
Lemma~D. It follows that, in the case of a finite spectrum, there must be a
value $\l\in Sp(T)$ that is infinitely singular for $T$.

We have therefore proved the following theorem.

\proclaim Theorem. If $X$ is a complex H.I. Banach space, then every bounded
linear operator $T$ from $X$ into $X$ can be written $T=\l I+S$, where
$\l\in\CC$ and $S$ is strictly singular. The spectrum of $T$ is finite, or
consists of $\l$ and a sequence $(\l_n)$ of eigenvalues with finite
multiplicity converging to $\l$.

\proclaim Corollary. A complex H.I. space $X$ is not isomorphic to any proper
subspace, and in particular is not isomorphic to its hyperplanes.

\bigbreak

\noindent{\bf \S 5. Further properties.}
\medskip

We shall now show how to pass from the complex case back to the real case.
The following lemma will be useful; it was shown to us by R.~Haydon.

\proclaim Lemma. Suppose $X$ is a real HI-space and
$T$ a bounded linear operator from $X$ into itself. If we denote by
$S$ the natural extension of $T$ to the complexification of $X$,
then the spectrum of $S$ is invariant by conjugation, and the part
in the upper complex plane is finite or consists of a converging
sequence.

As before, this lemma implies that there exists no isomorphism from $X$
onto a proper subspace.
\medskip

\Proof  If $\lambda\notin F_S$ is real,
there exists for every $\epsilon>0$ a
(real) infinite dimensional subspace $Y_\epsilon$ of $X$ such that
$\|T-\lambda Id_{Y_\epsilon}\| < \epsilon$
on $Y_\epsilon$. Since $X$
is HI, it follows that $\CC\setminus F_S$ contains at most one
real element. Let now $\lambda, \mu\notin F_S$, and
$\mu\notin\{\lambda, \bar\lambda\}$. We may assume that $\lambda$
is not real. Let
$$ T_\lambda = T^2 - 2Re \lambda T + |\lambda|^2 Id.$$
Then $(S-\bar\lambda Id)(S-\lambda Id)(x + iy) =
T_\lambda x + i \, T_\lambda y$; for every $\epsilon>0$
it is thus possible to find an infinite dimensional subspace
$Y_\epsilon$ of $X$ such that $\|T_\lambda \| < \epsilon$
on $Y_\epsilon$. Since $X$ is HI, we may assume the same
for $T_\mu$ on the same $Y_\epsilon$.
Now, $T_\lambda - T_\mu = aT + bId$
for some $a,b\in\R$, not both $0$, and it
has norm less than $2\epsilon$ on $Y_\epsilon$.
Thus $a\ne 0$. We obtain  that $T$ is nearly equal to
$(-b/a) Id$ on $Y_\e$; since $T_\l$ is nearly $0$
on $Y_\e$, we get easily
that $-b/a$ must be a root of
the polynomial $(X-\bar\lambda)(X-\lambda)$,
which is of course impossible.

  We know therefore that $\CC\setminus F_S$ contains at most
a pair $(\lambda, \bar\lambda)$, and the rest of the proof
is as in section~IV.

\medskip

\proclaim Theorem V.1. If $X$ is a real HI space (for example,
$X$ could be the real version of the space from section~III),
then $X$ is not isomorphic to any proper subspace. In particular,
$X$ is not isomorphic to its hyperplanes.

\bigbreak

\centerline{\cbf References}
\medskip

{\ninepoint

\myitem{[BP]} C. Bessaga and A.~Pe\pl czy\'nski,
{\it A generalization of results of R.C. James concerning
absolute bases in Banach spaces},
Studia Math. {\bf 17} (1958), 165--174.

\myitem{[CS]} P. Casazza, T. Shura,
{\it Tsirelson's space},
Lecture Notes in Math. {\bf 1363}, Springer Verlag.

\myitem{[D]} M. Day,
{\it Normed linear spaces},
Springer Verlag.

\myitem{[E]} P. Enflo, 
{\it A counterexample to the approximation property in
Banach spaces}, 
Acta Math. {\bf 130} (1973), 309--317. 

\myitem{[G]} T. Gowers, 
{\it A solution to the Banach hyperplane problem},  

\myitem{[LT]} J. Lindenstrauss and L. Tzafriri, 
``Classical Banach Spaces, I,'' 
Springer Verlag, 1977. 

\myitem{[MR]} B. Maurey and H. Rosenthal, 
{\it Normalized weakly null sequences with no unconditional 
subsequence}, 
Studia Math. {\bf61} (1977), 77--98. 

\myitem{[OS]} E. Odell and T. Schlumprecht,
{\it The distortion problem},

\myitem{[S1]} T. Schlumprecht, 
{\it An arbitrarily distortable Banach space},

\myitem{[S2]} T. Schlumprecht, 
{\it A complementably minimal Banach space not containing $c_0$ or $\ell_p$}, 

\myitem{[T]}  B.S. Tsirelson,
{\it Not every Banach space contains $\ell_p$ or $c_0$}, 
Funct. Anal. Appl. {\bf 8} (1974), 138--141. 

}

\bye